\documentclass[oneside,11pt]{article} %,draft,openright]

\usepackage{epsfig,graphicx}
\usepackage{latexsym,amssymb}
\usepackage{setspace,cite} %double spacing for text, single for captions,
\usepackage{amsmath, amssymb, amsthm}
\usepackage{graphicx,color}

\usepackage{fancyhdr}

\usepackage{a4}
\usepackage{bbm}

\newcommand{\B}[1]{\mathbb{#1}}
\usepackage{amsmath}
\usepackage{amssymb}
\usepackage{euscript}
\usepackage{amsthm}
\usepackage[all]{xy}
\newtheorem{prop}{Proposition}[section]

\newtheorem{cor}[prop]{Corollary}
\newtheorem{lemma}[prop]{Lemma}
\newtheorem{thm}[prop]{Theorem}

\theoremstyle{definition}
\newtheorem{defn}[prop]{Definition}

\theoremstyle{definition}

\theoremstyle{definition}
\newtheorem*{note}{Note}
\theoremstyle{definition}

\title{Line Bundles over quantum tori} 
\author{L. Taylor}

\begin{document}
\maketitle

\begin{abstract}
Quantum tori are a class of non-Hausdorff spaces proposed to take a central role in Manin's proposed theory of Real Multiplication, which provides a framework in which to study Hilbert's twelfth problem for real quadratic fields.  We define a notion of line bundles on such spaces, and prove a structure theorem concerning isomorphism classes of line bundles analogous to the Appel-Humbert Theorem for Complex Tori.  We introduce the notion of the Chern class of a line bundle and prove its relation to the Heisenberg group associated to such an object. Finally we identify the theta functions associated to these line bundles and consider their context within a solution to the generation of class fields over real quadratic fields. Finally we remark on the context of line bundles within Model Theory.
\end{abstract}

\section{Introduction}\label{intro3}

In \cite{Manin}, Manin proposed a framework in which to approach Hilbert's twelfth problem for real quadratic fields - an explicit class field theory for such fields.  An explit class field theory has been achieved for imaginary quadratic fields using the theory of elliptic curves with Complex Multiplication.  In Manin's theory of ``Real Multiplication'', a family of topological spaces known as \emph{quantum tori} are proposed to play a fundamental role.   

\begin{defn}
Let $L$ be a dense subgroup of $\B{R}$ of rank two.  We say that $L$ is a \emph{pseudolattice}.  The quantum torus associated to $L$ is defined to be the topological space $\B{R}/L$.
\end{defn}

A \emph{complex torus} is defined to be a topological space $\mathbb{C}/\Lambda$ for some lattice $\Lambda$.
Central to the theory of Complex Multiplication is the existence of meromorphic elliptic functions\footnote{An elliptic function for $\Lambda$ is a function $f$ on $\B{C}$ such that $f(v+\lambda)=f(v)$ for all $v \in \B{C}$, $\lambda \in \Lambda$.} for $\Lambda$, such as the \emph{Weierstrass $\wp$-function} that provides an equivalence between the categories of complex tori and elliptic curves over $\mathbb{C}$. \\ 

An analogous function for a quantum torus $Z_L$ would be a meromorphic function which was periodic with respect to the pseudolattice $L$.  It is clear that any such function $f$ is constant:  If $f$ does not have a pole at $v \in \B{C}$, then there exists a neighbourhood $U$ of $v$ containing a sequence on which $f$ is constant which has an accumulation point in $U$.  Hence $f$ is constant on $U$, and therefore is constant on $\mathbb{C}$.  \\

Elliptic functions can be viewed as quotients of \emph{theta functions} on the complex torus - functions satisfying certain periodicity conditions with respect to the lattice $\Lambda$.  These functions can be viewed as sections of line bundles over complex tori, and certain classes of these functions characterise line bundles up to isomorphism \cite{Swan}.  This was the original motivation behind the results in this paper. 

\subsection{Statement of results and overview}

We define the notion of a holomorphic line bundle $\mathcal{L}$ over a quantum torus $Z_L$, drawing analogies from the theory of such objects over complex tori.  Our cohomological definition admits a natural definition of isomorphism between line bundles, and identifies isomorhpism classes of line bundles over $Z_L$ as elements of a certain cohomology group $H^1(L,\mathcal{H}^\ast)$.\\

We note here that line bundles over quantum tori have been defined and studied by Polishchuck \cite{PolishchukI,PolishchukII}, Schwarz and Astashkevich \cite{SchwarzI,SchwarzIII} from the perspective of Noncommutative Geometry.  It would be interesting to see how the notion developed in this paper is related to that of the latter papers.\\

In \S\ref{AH} we prove a structure result concerning the group of isomorphism classes of line bundles over a quantum torus $Z_L$.  The proof of this result introduces the notion of the Chern class $Ch(\mathcal{L})$ of a line bundle, which is defined to be the image of $\mathcal{L}$ under a map (to be defined)
\begin{equation}\label{chmap}
Ch: H^1(L,\mathcal{H}^\ast) \rightarrow \textrm{Alt}^2(L,\B{Z}), 
\end{equation} 
where $\textrm{Alt}^2(L,\B{Z})$ denotes the group of integral valued alternating forms on $L$.
\begin{thm}\label{theorysplit}
Let $Z_L$ be a quantum torus.  Then we have an isomorphism
$$H^1(L,\mathcal{H}^\ast) \cong \B{C}^\times \oplus \textrm{Alt}^2(L,\B{Z}).$$
\end{thm}

In \S\ref{AHQT} we use the above decomposition of $H^1(L,\mathcal{H}^\ast)$ to prove a result analogous to the Appel-Humbert Theorem for complex tori, which classifies isomorphism classes of line bundles in terms of a hermitian form and a type of character:

\begin{thm}\label{Apple}
Let $L$ be a pseudolattice.  Let $P(L)$ denote the group whose elements are pairs $(\chi,E)$ such that $E \in \textrm{Alt}^2(L,\B{Z})$ and $\chi:L \rightarrow \B{C}^\times$ is such that for all $l_1,l_2 \in L$
$$\chi(l_1+l_2)=\chi(l_1)\chi(l_2)e^{2 \pi i E(l_1,l_2)}.$$ 
The group law in $P(L)$ is given by
$$(\chi_1,E_1)(\chi_2,E_2)=(\chi_1\chi_2,E_1+E_2).$$
Then we have an ismorphism $H^1(L,\mathcal{H}^\ast) \cong P(L)$.
\end{thm}

In \S\ref{Heisenberg} we define the \emph{Heisenberg group} associated to a line bundle $\mathcal{L}$ over a quantum torus.  The development of this group gives rise to a certain subgroup $K(\mathcal{L})$ associated to each line bundle $\mathcal{L}$, together with an alternating pairing
$$e^\mathcal{L}: K(\mathcal{L}) \times K(\mathcal{L}) \longrightarrow \B{C}^\times.$$
We prove a result which shows how these notions are closely linked to the Chern class of a line bunde:

\begin{thm}\label{dichotomy}
Let $\mathcal{L}$ be a line bundle over a quantum torus $Z_L$.  
\begin{enumerate}
\item \label{nontriviald}If $Ch(\mathcal{L})$ is nontrivial, then $K(\mathcal{L})$ is finite.  In this case we have $$Ch(\mathcal{L})(a\omega_1+b\omega_2,c\omega_1+d\omega_2)=s_\eta(ad-bc)$$
for some $s_\eta \in \B{Z}$, and the following isomorphism holds:
$$K(\mathcal{L}) \cong (\B{Z}/s_\eta\B{Z})\times (\B{Z}/s_\eta\B{Z}).$$
\item \label{triviald}The following statements are equivalent:
\begin{enumerate}
\item \label{p}$Ch(\mathcal{L}) = 0$;
\item \label{q}$K(\mathcal{L})=Z_L$;
\item \label{r}$e^\mathcal{L} \equiv 1$.
\end{enumerate}
\end{enumerate}
\end{thm}

In the final section we discuss the existence of \emph{theta} functions corresponding to the line bundles whose theory we have developed.  We show that our notion of \emph{holomorphic} line bundles may be too strict for possible applications to Hilbert's problem for real quadratic fields, but observe from the work of Shintani \cite{ShintaniIII} that a weakening of this condition gives rise to functions that can be used to generate abelian extensions of real quadratic fields.

\section{Defining Line Bundles over quantum tori}\label{defining}

Let $L \subseteq \B{R}$ be a pseudolattice, and let $Z_L$ denote the associated quantum torus.  Our aim is to define the notion of a line bundle over $Z_L$.  The classical definition is given below:
\begin{defn}\label{line}
Let $X$ be a topological space.  A complex line bundle $\mathcal{L}$ over $X$ is a topological space $\mathcal{L}$ equipped with a projection 
$$\Pi: \mathcal{L} \rightarrow X$$
such that
\begin{itemize}
\item For each $x \in X$, $\Pi^{-1}(x)$ is a one dimensional $\B{C}$-vector space;
\item For each $x \in X$ there exists an open neighbourhood $U$ of $x$ such that $\Pi^{-1}(U) \cong U \times \B{C}$.
\end{itemize}
\end{defn}
Due to the non-Hausdorff nature of quantum tori, the only such line bundle over a quantum torus $Z_L$ is the trivial line bundle $\mathcal{L} \cong Z_L \times \mathbb{C}$.\\

Let $X_\Lambda$ be the complex torus defined by $\B{C}/\Lambda$ for some lattice $\Lambda$.  A description of holomorphic line bundles over $X_\Lambda$ exists as the group of cocycles $Z^1(\Lambda, \mathcal{H}^\ast)$, where $\mathcal{H}$ denotes the ring holomorphic functions on $\B{C}$.  Isomorphism classes of such line bundles are characterised by the cohomology group $H^1(\Lambda, \mathcal{H}^\ast)$.\\

This motivates the following definition:

\begin{defn}
Let $Z_L$ be a quantum torus.  A holomorphic line bundle $\mathcal{L}$ over $Z_L$ is an element of $Z^1(L,\mathcal{H}^\ast)$.  We say that two holomorphic line bundles are isomorphic if they have the same image in $H^1(L,\mathcal{H}^\ast)$, and denote by $[\mathcal{L}]$ the isomorphism class of $\mathcal{L}$.  We denote the law of composition both in $Z^1(L,\mathcal{H}^\ast)$ and $H^1(L,\mathcal{H}^\ast)$ by $\otimes$.  
\end{defn}
\begin{note}
In this paper we will refer to holomorhic line bundles simply as line bundles.  We will keep using the classical notation $\mathcal{L}$ for a line bundle over a quantum torus $Z_L$, and when we refer explicitly to the coycle we will denote this by $A_l^\mathcal{L}(v)$.  Our motivation for this is two fold:
\begin{itemize}
\item By using this notation it is easy to compare our results with the coresponding results for complex tori;
\item We wish to maintain the philosophy that line bundles are topological objects -  the line bundle $\mathcal{L}$ can be thought of as the topological space $\mathcal{L}:=(\mathbb{R} \times \mathbb{C})/L$, where the action of $L$ is given by 
$$
\begin{array}{rcl}
L \times \mathbb{R} \times \mathbb{C} &\rightarrow& \mathbb{R} \times \mathbb{C}\\
(l,v,z) & \mapsto & (v+l,A_l^\mathcal{L}(v)z).
\end{array}$$
\end{itemize}

\end{note}
\begin{prop}\label{existence}
Let $Z_L$ be a quantum torus.  Then there exist nontrivial line bundles on $Z_L$.
\end{prop}
\begin{proof}
Suppose $L=\B{Z}\omega_1+\B{Z}\omega_2$, and let $\mathcal{L}$ be the line bundle defined by the cocycle
$$
\begin{array}{rcl}
 L  \times \B{C} &\longrightarrow& \B{C}^\ast\\
(l,v) & \mapsto & A_l^\mathcal{L}(v) = e^{-\frac{\pi i}{\omega_1}[ b^2 {\omega_2} - 2  b v]},
\end{array}
$$
where $l=a\omega_1+b\omega_2$.  \\

Suppose  $A_l^\mathcal{L}(v)$ is cohomologically trivial.  Then there exists
$h \in \mathcal{H}^\ast$ such that
$$A_l^\mathcal{L}(v)=\frac{h(v+l)}{h(v)}$$
for all $v \in \B{C}, l \in L$.
Since $h$ is nonvanishing we may write $h(v)=e^{\pi i g(v)}$ for some holomorphic function
$g(v)$ which satisfies the following periodicity relations
\begin{eqnarray*}
g(v+\omega_1)-g(v)&=&2m\\
g(v+\omega_2)-g(v)&=&-\frac{\omega_2+2v}{\omega_1}
\end{eqnarray*}
for some $m \in \B{Z}$.  The holomorphic function $k(v):=g(v)-2mv/\omega_1$ satisfies the following periodicity conditions:
\begin{eqnarray}
k(v+\omega_1)-k(v)&=&0 \label{period1}\\
k(v+\omega_2)-k(v)&=&-\frac{(2m+1)\omega_2+2v}{\omega_1}\label{period2}
\end{eqnarray}
Consider the continued fraction expansion of $\theta:=\omega_2/\omega_1$, and let
$p_n/q_n$ be the convergents (see \cite{Burton}).  Then the sequence $q_n \rightarrow \infty$, but
\begin{equation}\label{conv5}
\left|p_n\omega_1-q_n \omega_2 \right|<\frac{\left| \omega_1 \right| }{  q_n} \rightarrow 0
\end{equation}
as $n \rightarrow \infty$.  Consider the sequence $x_n=p_n \omega_1-q_n \omega_2$. By the continuity of $k$, and \eqref{conv5} we have $k(x_n) \rightarrow k(0)$.  Hence by \eqref{period2} we obtain

\begin{equation}
\begin{array}{rcl}
0 &=& \lim_{n \rightarrow \infty}\left|k(x_n)-k(0) \right|\\
&=& \lim_{n \rightarrow \infty}\left|k(q_n \omega_2)-k(0) \right| \\
&=& \lim_{n \rightarrow \infty} q_n\left|\theta\right| \left| (2m+1)+(q_n+1)\right|.\label{infinite}
\end{array}
\end{equation}
The expression in the final line of \eqref{infinite} tends to $\infty$ as $n \rightarrow \infty$, which is a contradiction.\newline

Therefore $A_l^\mathcal{L}(v)$ is a nontrivial element of $H^1(L,\mathcal{H}^\ast)$.
\end{proof}

\section{An Appel-Humbert Theorem for quantum tori} \label{AH}

The goal of this section is to prove Theorem \ref{theorysplit}. The proof relies on the definiton of the Chern class of a line bundle.

\subsection{The Chern Class of a Line Bundle}

Let $L$ be a pseudolattice, and let $\mathcal{H}$ denote the ring of holomorphic functions on $\B{C}$.  We have a natural action of $L$ on $\mathcal{H}$ given by 
\begin{equation}\label{Laction}
l.f(v)=f(v+l).
\end{equation}  
This equation also defines an action of $L$ on $\mathcal{H}^\ast$, and we have the following short exact sequence of $L$-modules
$$0 \rightarrow \B{Z} \rightarrow \mathcal{H} \rightarrow \mathcal{H}^\ast \rightarrow 0,$$
where the action on $\B{Z}$ is trivial.  Taking the connecting map in the long exact sequence obtained from this we obtain a homomorphism
$$\partial: H^1(L,\mathcal{H}^\ast) \longrightarrow H^2(L,\B{Z}).$$
\begin{defn}
Let $\mathcal{L}$ be a line bundle.  The Chern class of $\mathcal{L}$ is given by $\partial(\mathcal{L}) \in H^2(L,\mathcal{H}^\ast)$.
\end{defn}
The map $\partial$ can be defined on cocycles, and shown to map coboundaries to coboundaries.  We let $\hat \partial$ denote the map on cocycles which induces the map $\partial$ on cohomology groups.  Given a line bundle $\mathcal{L}$ represented by a cocycle $A_l^\mathcal{L}(v) \in Z^1(L,\mathcal{H}^\ast)$, the cohomology theory provides us with an explicit formula for the image of $\mathcal{L}$ under $\hat \partial$. \\

If $A_l^\mathcal{L}(v)=e^{2 \pi i a_\mathcal{L}(l,v)}$, then $\hat \partial(A)$ is a function on $L \times L$ taking values in $\B{Z}$ and has the explicit formula (see the proof of Theorem 2.1.2 in \cite{Lange})
$$\hat \partial (\mathcal{L})(l_1,l_2)=a_\mathcal{L}(l_1+l_1,v)-a_\mathcal{L}(l_1,v)-a_\mathcal{L}(l_2,v+l_1)$$
for some $v \in \B{C}$.  This is well defined since the cocycle condition satisfied by $A_l^\mathcal{L}(v)$ implies that this expression is independent of the choice of $v$.\newline

Define a map
$$\begin{array}{rcl}
\alpha: Z^2(L,\B{Z}) & \rightarrow & \textrm{Alt}^2(L,\B{Z})\\
P & \mapsto & \alpha(P), 
\end{array}
$$
where $\alpha(P)$ is defined by 
$$\alpha(P)(\omega_1,\omega_2) = P(\omega_1,\omega_2)-P(\omega_2,\omega_1).$$
This induces a well defined map (also denoted by $\alpha$) from $H^2(L,\B{Z})$ to the space $\textrm{Alt}^2(L,\B{Z})$ of alternating forms on $L$.
\begin{lemma} \label{alt}
The map $\alpha: H^2(L,\B{Z}) \rightarrow \textrm{Alt}^2(L,\B{Z})$ is an isomorphism.
\end{lemma}
\begin{proof}
See Lemma 2.1.3 of \cite{Lange}.
\end{proof}

The composition of $\hat \partial$ with $\alpha$ establishes a homomorphism
\begin{equation}\label{Ch}
\begin{array}{rcl}
Ch:H^1(L,\mathcal{H}^\ast) &\rightarrow &\textrm{Alt}^2(L,\B{Z})\\
~[\mathcal{L}] & \mapsto & \alpha \circ  [\hat\partial {(\mathcal{L})}]
\end{array}
\end{equation}  
An explicit expression for this map
is given by (Theorem 2.1.2 of \cite{Lange})
\begin{equation} \label{chernequation}
Ch(\mathcal{L})(l_1,l_2)=a_\mathcal{L}(l_2,v+l_1)+a_\mathcal{L}(l_1,v)-a_\mathcal{L}(l_2,v)-a_\mathcal{L}(l_1,v+l_2).
\end{equation}
Analysing this map in more detail will enable us to determine the structure of $H^1(L,\mathcal{H}^\ast)$.

\subsection{Surjectivity of ${Ch}: H^1(L,\mathcal{H}^\ast) \rightarrow \textrm{Alt}^2(L,\B{Z})$}
\label{chsur}

In this section we study the image of the homomorphism $Ch$ defined in \eqref{Ch}.  We will prove the following result:

\begin{prop}\label{sursplit}
The map $Ch:H^1(L,\mathcal{H}^\ast) \rightarrow \textrm{Alt}^2(L,\B{Z})$ is surjective.  Furthermore, there exists a map 
$\sigma: \textrm{Alt}^2(L,\B{Z}) \rightarrow H^1(L,\mathcal{H}^\ast)$ such that $Ch \circ \sigma$ is the identity on $\textrm{Alt}^2(L,\B{Z})$.
\end{prop}
\begin{proof}
The proof relies on the construction of elements of $Z^1(L,\mathcal{H}^\ast)$ which are similar to that introduced in the proof of Proposition \ref{existence}.\newline

If $L=\mathbb{Z}\omega_1+\mathbb{Z}\omega_2$ then every $\eta \in \textrm{Alt}^2(L,\B{Z})$ is of the form
$$\eta(a\omega_1+b\omega_2,c\omega_1+d\omega_2)=s_\eta(ad-bc)$$
for some $s_\eta \in \B{Z}$.  The assignment $\eta \mapsto s_\eta$ gives a bijection between $\textrm{Alt}^2(L,\B{Z})$ and $\B{Z}$.\newline

For each $\eta \in \textrm{Alt}^2(L,\B{Z})$ define $\mathcal{L}_\eta$ to be the line bundle associated to the coycle
\begin{equation}\label{sigmasp}
\hat \sigma(\eta)_l(v):=e^{{s_\eta}\frac{\pi i}{\omega_1}[b^2 \omega_2+2bv]}
\end{equation}
where $l=a\omega_1+b\omega_2$.  We define $\sigma(\eta):=[\mathcal{L}_\eta]$.\\

By the definition above we have
$\hat \sigma(\eta)_l(v)=e^{2 \pi i \Sigma_\eta(l,v)}$ where
$$\Sigma_\eta(l,v)=s_\eta \frac{1}{2 \omega_1}[b^2 \omega_2+2bv].$$
We calculate $Ch(\sigma(\eta))$ using \eqref{chernequation}.  Let $l_1=a\omega_1+b\omega_2$, $l_2=
c\omega_1+d\omega_2 \in L$.  Then
$$
\begin{array}{rcl}
Ch(\sigma(\eta))(l_1,l_2) & = & \Sigma_\eta(l_2,v+l_1)+\Sigma_\eta(l_1,v)-\Sigma_\eta(l_2,v)
-\Sigma_\eta(l_1,v+l_2)\\
&=&\frac{s_\eta}{2 \omega_1}[d^2 \omega_2+2d(v+a\omega_1+b\omega_2)-d^2 \omega_2-2dv]\\
& & \qquad \qquad -\frac{s_\eta}{2 \omega_1}[b^2 \omega_2+2b(v+c\omega_1+d\omega_2)-b^2 \omega_2-2bv]\\
&=& s_\eta(ad-bc)\\
&=& \eta(l_1,l_2).
\end{array}
$$

\end{proof}

\subsection{The kernel of ${Ch : H^1(L,\mathcal{H}^\ast) \rightarrow \textrm{Alt}^2(L,\B{Z})}$.}\label{ker}

Traditionally the kernel of the Chern map is denoted by the following:
\begin{defn}
Let $\textrm{Pic}^0(Z_L)=\{[\mathcal{L}] \in H^1(L,\mathcal{H}^\ast): Ch([\mathcal{L}])=0 \}$.
\end{defn}

We now show that every line bundle with trivial Chern class over a quantum torus is isomorphic to one represented by a \emph{constant} cocycle.  

\begin{prop} \label{constant}
Let $\mathcal{L}$ be a line bundle over a quantum torus such that $Ch([\mathcal{L}])=0$.  Then there exists a trivial line bundle $\mathcal{L}'$ such that $\mathcal{L} \otimes \mathcal{L}'$ is given by a cocycle that is constant. 
\end{prop}
\begin{proof}
The proof of this result involves unravelling the definition of the map $Ch$.
Let $A_l^\mathcal{L}(v)=e^{2 \pi i a_\mathcal{L}(l,v)}$ for some function $a: L \times \B{C} \rightarrow \B{C}$
holomorphic for fixed $l \in L$.  Since $A_l^\mathcal{L}(v)$ satisfies the cocycle condition we have
\begin{equation}\label{cocycle}
a_\mathcal{L}(l_1+l_2,v)-a_\mathcal{L}(l_1,v+l_2)-a_\mathcal{L}(l_2, v) \equiv 0 \pmod{\B{Z}}.
\end{equation}
The image of $A_l^\mathcal{L}(v)$ under the map $\hat \partial: Z^1(L,\mathcal{H}^\ast) \rightarrow
Z^2(L,\B{Z})$ is trivial. By \eqref{chernequation} this implies that
\begin{equation}\label{chern}
a_\mathcal{L}(l_2,v+l_1)+a_\mathcal{L}(l_1, v)-a_\mathcal{L}(l_2,v)-a_\mathcal{L}(l_1,v+l_2)=0
\end{equation}
for every $v \in \B{C}$ and for all $l_1,l_2 \in L$.  Define $h(v):=a_\mathcal{L}(0, v)$.  Then
\begin{equation}\label{bun}
\begin{array}{rcl}
a_\mathcal{L}(l, v)+h(v+l)-h(v) & = & a_\mathcal{L}(l, v)+a_\mathcal{L}(0, v+l)-a_\mathcal{L}(0, v)\\
& \equiv & a_\mathcal{L}(l, v)-a_\mathcal{L}(0, v) \pmod{\B{Z}}\\
& \equiv & a_\mathcal{L}(l, 0)-a_\mathcal{L}(0,  l) \pmod{\B{Z}}.
\end{array}
\end{equation}
The second line follows from \eqref{cocycle} by putting $l_1=0$, $l_2=l$.  Using these
same substitutions in \eqref{chern}, together with $v=0$ yields the final line.  This is valid since \eqref{chern} is independent of $v$.  Let $H(v)=e^{2 \pi i h(v)}$, and let $\mathcal{L}'$ be the trivial line bundle represented by the cocycle $B_l(v):=H(v+l)H(v)^{-1}$.  Then $\mathcal{L} \otimes \mathcal{L}'$ is represented by the coycle $K_l(v):=A_l^\mathcal{L}(v)B_l(v)$ is independent of $v$.
\end{proof}

The previous result implies the existence of a homomorphism 
\begin{equation}\label{homc}
C: \{\mathcal{L} \in Z^1(L,\mathcal{H}^\ast) : [\mathcal{L}] \in \textrm{Pic}^0(Z_L)\} \rightarrow \textrm{Hom}(L,\B{C}^\times).
\end{equation}
By the proof of Proposition \ref{constant}, $C$ is constant on the cohomology classes, and so can be viewed as a function on $\textrm{Pic}^0(Z_L)$.\\

Note that any element of $\textrm{Hom}(L,\B{C}^\times)$ can be viewed as an element of $Z^1(L,\mathcal{H}^\ast)$.

\begin{lemma}\label{prep}
Let $\phi \in \textrm{Hom}(L,\B{C}^\times)$.  Then viewed as an element of $Z^1(L,\mathcal{H}^\ast)$, $\phi$ is cohomologically trivial if and only if there exists $k \in \B{C}^\times$ such that $\phi(l)=k^l$ for all $l \in L$.
\end{lemma}
\begin{proof}
If $\phi$ is cohomologically trivial, then there exists $h \in \mathcal{H}^\ast$ such that for all $l \in L$
\begin{equation}\label{phitriv}
\phi(l)=\frac{h(v+l)}{h(v)}.
\end{equation}
Since $h$ is continuous, $\phi$ extends to a continuous homomorphism $\hat \phi$ of $\B{R}$ in to $\B{C}$.  Then for any $r \in \B{R}$ we have $\hat \phi(r)=\hat \phi(1)^r:=e^{r \log(\hat \phi(1))}$ for some choice of branch of the logarithm function.  Hence for all $l \in L$, $\phi(l)=k^l$, where $k=\phi(1)$.\\

Conversely, if $\phi(l)=k^l$ for all $l \in L$, define $h(v)=k^v$.  Then $h$ satisfies \eqref{phitriv} and hence $\phi$ is cohomologically trivial. 
\end{proof}

\begin{prop}\label{pic0}
$\textrm{Pic}^0(Z_L) \cong \B{C}^\times$.
\end{prop}
\begin{proof}
Suppose $L=\B{Z}\omega_1+\B{Z}\omega_2$, and let $\mathcal{L}$ be such that $Ch(\mathcal{L})=0$.  Let $\phi = C(\mathcal{L}) \in \textrm{Hom}(L,\mathcal{H}^\ast)$, and define $k=\phi(\omega_1)^{-1}$.  Denote by $\hat \phi$ the element $\hat \phi(l)=\phi(l)k^l$ of $\textrm{Hom}(L,\B{C}^\times)$. By Lemma \ref{prep} $\phi$ and $\hat \phi$ are cohomologous and we have $\hat \phi(\omega_1)=1$.
Define the homomorphism 
\begin{equation}\label{CC}
\begin{array}{rcl}
\textrm{Pic}^0(Z_L) & \rightarrow & \B{C}^\times\\
\phi:=C([\mathcal{L}]) & \mapsto & \hat \phi(\omega_2).
\end{array}
\end{equation}
This is well-defined: If $\phi=C([\mathcal{L}])$ is cohomologically trivial then we have $\phi(l)=\phi(1)^l$ for all $l \in L$, and therefore $\hat \phi(l)=1$ for all $l \in L$.\\

It is clear that this map is injective and surjectivity also follows: for any $c \in \B{C}^\times$ define a homomorphism $\psi$ of $L$ by $\psi(a\omega_1+b\omega_2)=c^b$.  Then the image of $\psi$ is $c$.

\end{proof}

\begin{proof}[of Theorem \ref{theorysplit}]
By Proposition \ref{sursplit} we have the following short split exact sequence.
$$\xymatrix{0 \ar[r] & \textrm{Pic}^0(Z_L) \ar[r] & H^1(L,\mathcal{H}^\ast) \ar[r] & \textrm{Alt}^2(L,\B{Z}) \ar[r] & 0}.$$
By Proposition \ref{pic0} we have an isomorphism $\textrm{Pic}^0(Z_L) \cong \B{C}^\times$.  The theorem follows from the theory of split exact sequences.
\end{proof}

\section{An Appell-Humbert Theorem for quantum tori}\label{AHQT}

We now formulate a classification of isomorhpism classes of line bundles over quantum tori analogous to the Appell-Humbert Theorem for complex tori.  \\

Let $\mathcal{L} \in Z^1(L,\mathcal{H}^\ast)$, and fix a $\B{Z}$-basis $\{\omega_1,\omega_2\}$ of $L$.  Then with the notation of Proposition \ref{sursplit}, we have 
$$[\mathcal{L}]^{-1} \otimes \sigma(Ch([\mathcal{L}])) \in \textrm{Pic}^0(Z_L).$$
Via the isomorphism of \eqref{CC} this defines an element $c_\mathcal{L} \in \B{C}^\times$.  By Theorem \ref{theorysplit}, the isomorphism class of $\mathcal{L}$ is characterised by the pair $(c_\mathcal{L},Ch([\mathcal{L}]))$.  We denote by $\gamma_\mathcal{L}$ the element of $\textrm{Hom}(L,\B{C}^\times)$ given by $\gamma_\mathcal{L}(a\omega_1+b\omega_2)=c_\mathcal{L}^b$.\\

Recalling the definitions in Theorem \ref{Apple} we note that the maps $\alpha: z \mapsto (z,0)$ and $\beta(z,E)=E$ define the following short exact sequence
$$\xymatrix{0 \ar[r]^\alpha & \B{C}^\times \ar[r] & P(L) \ar[r]^\beta & \textrm{Alt}^2(L,\B{Z}) \ar[r] & 0}.$$

\begin{prop}
There exists a homomorphism $\Phi: H^1(L,\mathcal{H}^\ast) \rightarrow P(L)$.
\end{prop}
\begin{proof}
Recall from Proposition \ref{sursplit} that we have an isomorphism 
$$
\begin{array}{rcl}
\textrm{Alt}^2(L,\B{Z}) &\cong &\B{Z}\\
\eta & \mapsto & s_\eta.
\end{array}
$$
For a line bundle $\mathcal{L}$, we define $s_\mathcal{L} \in \B{Z}$ to be the image of $Ch([\mathcal{L}])$ under the above isomorphism.  We then define $\chi_\mathcal{L}: L \rightarrow \B{C}^\times$ by
$$\chi_\mathcal{L}(a\omega_1+b\omega_2)=e^{\pi i s_\mathcal{L}ab}.$$
For $l_1,l_2 \in L$ this satisfies 
$$\chi_\mathcal{L}(l_1+l_2)=\chi_\mathcal{L}(l_1)\chi_\mathcal{L}(l_2)e^{\pi i Ch(\mathcal{L})(l_1,l_2)}.$$
Using the decompositon of $H^1(L,\mathcal{H}^\ast)$ as in Theorem \ref{theorysplit} we define the homomorphism
$$
\begin{array}{rcl}
\Phi: H^1(L,\mathcal{H}^\ast) & \rightarrow & P(L)\\
(c_\mathcal{L},Ch([\mathcal{L}]) & \mapsto & (\gamma_\mathcal{L} \chi_\mathcal{L}, Ch([\mathcal{L}])).
\end{array}$$
\end{proof}

\begin{proof}[of Theorem \ref{Apple}]
If we can show that the following diagram is commutative then we are done by the snake lemma:
$$
\xymatrix{
0 \ar[r] & \textrm{Pic}^0(Z_L) \ar[d] \ar[r]& H^1(L,\mathcal{H}^\ast) \ar[d] \ar[r]^{Ch} & \textrm{Alt}^2(L,\B{Z}) \ar[r] \ar@{=}[d]& 0\\
0 \ar[r] & \B{C}^\times \ar[r]^\alpha & P(L) \ar[r]^\beta & \textrm{Alt}^2(L,\B{Z}) \ar[r] & 0
}
$$
where the map $\textrm{Pic}^0(Z_L) \rightarrow \B{C}^\times$ is this isomorphism of \eqref{CC}. It is easy to show that this diagram is commutative.
\end{proof}

This can be compared to the \emph{Appel-Humbert Theorem} for complex tori.  The differences in the statements of these theorems for quantum and complex tori occur because for complex tori the induced topology from $\B{C}$ on the lattice is discrete, where as the induced euclidean topology from $\B{R}$ on to the pseudolattice is not.

\section{The Heisenberg Group}\label{Heisenberg}

Suppose $X$ is a topological space endowed with a group law $+$.  Given $x \in X$ we have a natural ``translation by $x$'' map
$$
\begin{array}{rcl}
T_x:X &\rightarrow& X\\
y & \mapsto & y+x.
\end{array}$$ 
Given a line bundle $\mathcal{L}$ over $X$ we let $T^\ast(\mathcal{L})$ denote the pull back of $\mathcal{L}$ with respect to this map.

\begin{defn}
Let $X$ be a topological group, and $\pi_{\mathcal{L}}: \mathcal{L} \rightarrow X$ a line bundle over $X$.  Define $$K(\mathcal{L}):=\{x \in X : T_x^\ast(\mathcal{L}) \cong \mathcal{L}\}.$$
\end{defn}

The group $K(\mathcal{L})$ is fundamental in defining the Heisenberg Group associated to a line bundle $\mathcal{L}$ over a topological space $X$. 

\begin{defn}[Heisenberg Group]
Let $\pi_\mathcal{L}: \mathcal{L} \rightarrow X$ be a line bundle over a topological group $X$.
The Heisenberg Group $H(\mathcal{L})$ of $\mathcal{L}$, is defined to be the set of pairs $(x,\phi)$ such that
\begin{itemize}
\item $x \in K(\mathcal{L})$;
\item $\phi: \mathcal{L} \rightarrow T_x^\ast(\mathcal{L})$ is an isomorphism.
\end{itemize}
The group law on $H(\mathcal{L})$ is given by
$$(x_1,\phi_1).(x_2,\phi_2)=(x_1+x_2,T_{x_1}^\ast(\phi_2) \circ \phi_1). $$
\end{defn}

\begin{note}
In \cite{ManinIII}, Manin defines a notion of the Heisenberg Group associated to a quantum torus.  His treatment also follows an analogy for complex tori, but this time focusses on viewing the Heisenberg group as having a representation over the space of holomorphic theta functions.  It would be an interesting study to explore the relationship between these two notions of Heisenberg groups for quantum tori.
\end{note}

\begin{defn}
Let $\mathcal{L}$ be a line bundle over a quantum torus $Z_L$ represented by $A_l(v) \in Z^1(L,\mathcal{H}^\ast)$, and suppose $x \in Z_L$.  Furthermore let $\tilde x$ be such that $\pi(\tilde x)=x$. Then the pull back of $\mathcal{L}$ with respect to the pair $(T_x, \tilde x)$ is the line bundle $T^\ast_{x, \tilde x}(\mathcal{L})$ represented by the cocycle $({T_{x,\tilde x}^\ast A})_l(v)=A_l(v+\tilde x)$.
\end{defn}

\begin{prop}
Let $\mathcal{L}$ be a line bundle over $Z_L$ and suppose $\pi(\tilde x_1)=\pi(\tilde x_2)=x \in Z_L$.  Then $T^\ast_{x,\tilde x_1}(\mathcal{L}) \cong T^\ast_{x,\tilde x_2}(\mathcal{L}).$ 
\end{prop}
\begin{proof}
We have $\tilde x_1-\tilde x_2 =l' \in L$, and hence
$$\frac{A_l(v+\tilde x_1)}{A_l(v+\tilde x_2)}=\frac{A_{l}(v+\tilde x_2+l')}{A_{l}(v+\tilde x_2)} =\frac{A_{l'}(v+\tilde x_2+l)}{A_{l'}(v+\tilde x_2)}\in
B^1(L,\mathcal{H}^\ast).$$
\end{proof}

\begin{defn}
For a line bundle $\mathcal{L}$ over $Z_L$, define $[T^\ast_x(\mathcal{L})]$ to be the isomorphism class of $T^\ast_{x,\tilde x}(\mathcal{L})$ for some choice of $\tilde x$.
\end{defn}

\begin{prop}\label{heis}
Let $\mathcal{L}$ be a line bundle over a quantum torus $Z_L$ represented by a cocycle $A_l(v)$.  The Heisenberg group $H(\mathcal{L})$ is given by pairs $(x,h) \in Z_L \times \mathcal{H}^\ast$ such that for some $\tilde x$ with $\pi(\tilde x)=x$ we have
$$\frac{A_l(v+\tilde x)}{A_l(v)}=\frac{h(v+l)}{h(v)}.$$
The law of composition is given by 
$$(x_1,h_1(v)).(x_2, h_2(v))=(x_1+x_2,h_2(v+\tilde x_1)h_1(v)). $$
\end{prop}
\begin{proof}
This follows from the definitions.
\end{proof}

\begin{lemma} \label{uniqueho}
Given $\tilde x \in \mathbb{R}$, any function $h_{\tilde x}$ such that
$$\frac{A_l(v+\tilde x) }{ A_l(v)}=\frac{h_{\tilde x}(v+l) }{ h_{\tilde x}(v)}$$
for all $l \in L$ is unique up to multiplication by an element of $\B{C}^\times$.
\end{lemma}
\begin{proof} Suppose there were two such functions, $h_{1,\tilde x}$ and $h_{2,\tilde x}$.  Then 
$$\frac{h_{1,\tilde x}(v+l)}{h_{2,\tilde x}(v+l)}=\frac{h_{1,\tilde x}(v)}{h_{2,\tilde x}(v)} $$
for all $l \in L$.  The right hand side is independant of $l$, so the left hand side is.  Since $L$ is dense in $\B{R}$, and both sides are holomorphic, the above expression is constant.  
\end{proof}

\begin{prop}
Let $\mathcal{L}$ be a line bundle over a quantum torus $Z_L$.  Then we have an exact sequence
$$\xymatrix{1 \ar[r] & \B{C}^\times \ar[r]^{\iota \ \ } & H(\mathcal{L}) \ar[r]^p & K(\mathcal{L}) \ar[r] & 0.} $$
\end{prop}
\begin{proof}
Define $\iota: \B{C}^\times \rightarrow H(\mathcal{L})$ by $\iota(k)=(0,k) \in K(\mathcal{L}) \times \mathcal{H}^\ast$.  Lemma \ref{uniqueho} shows that this is injective.  The map $p:K(\mathcal{L}) \times \mathcal{H}^\ast \rightarrow K(\mathcal{L})$ defined by $p(x,h)=x$ is clearly surjective.  It is also clear that $p \circ \iota (k)=0$ for all $k \in \B{C}^\times$.
\end{proof}

This last result establishes that $H(\mathcal{L})$ is a \emph{central extension} of $K(\mathcal{L})$.  According to the theory of such extensions we have an alternating pairing $K(\mathcal{L})$ defined by:

$$
\begin{array}{rcl}
e^\mathcal{L}: K(\mathcal{L}) \times K(\mathcal{L}) & \rightarrow & \B{C}^\times \\
(x_1,x_2) &\mapsto &\iota^{-1}(g_1g_2g_1^{-1}g_2^{-1}),
\end{array}$$
where $g_i=(x_i,h_i) \in H(\mathcal{L})$ for $i=1,2$.

\subsection{An alternating pairing on $\Lambda(\mathcal{L})$}

For every line bundle $\mathcal{L}$ over the quantum torus $Z_L$, we have the following objects:
\begin{itemize}
\item A subgroup $K(\mathcal{L})$ of $Z_L$;
\item An alternating pairing $e^\mathcal{L}$ on $K(\mathcal{L})$;
\item An alternating pairing $Ch(\mathcal{L})$ on $L$ taking integral values.
\end{itemize}

Our aim is to describe how these objects are related. \\

We make the following definition:

\begin{defn}
$\Lambda(\mathcal{L})=\{ \tilde x \in \B{R}: \pi(\tilde x) \in K(\mathcal{L})\}.$
\end{defn}

We will define an alternating pairing on $\Lambda(\mathcal{L})$ and show how this enables us to compute $e^\mathcal{L}$.

\begin{defn}\label{lambdadefine}
Fix $v \in \B{C}$ and a line bundle $\mathcal{L}$ over a quantum torus $Z_L$.  Let $x \in K(\mathcal{L})$ and fix $\tilde x \in \Lambda(\mathcal{L})$ such that $\pi(\tilde x)=x$.  By Lemma \ref{uniqueho} there exists a unique function $h_{\tilde x} \in \mathcal{H}^\ast$ such that $h_{\tilde x}(0)=1$ and
$$\frac{A_l(v+\tilde x)}{ A_l(v)}=\frac{h_{\tilde x}(v+l)}{ h(v)}.$$
Define a pairing on $\Lambda(\mathcal{L})$ (depending on $v$) by the following formula:
$$
\begin{array}{rcl}
H_v( \ , \ ): \Lambda(\mathcal{L}) \times \Lambda(\mathcal{L}) & \rightarrow & \B{C}^\ast\\
(\tilde x_1, \tilde x_2) & \mapsto & \frac{h_{\tilde x_2}(v+\tilde x_1)}{h_{\tilde x_2}(v)}.
\end{array}
$$

\end{defn}

\begin{prop} \label{c2}
For $\tilde x_1,\tilde x_2 \in \Lambda(\mathcal{L})$, and any $v \in \B{C}$ we have
$$e^\mathcal{L}(x_1, x_2)=H_v(\tilde x_1,\tilde x_2)H_v(\tilde x_2,\tilde x_1)^{-1}$$
where $x_i=\pi(\tilde x_i)$ for $i=1,2$.
\end{prop}
\begin{proof}
Firstly note that by the definition of the group law we have
$$(\tilde x,h_{\tilde x}(v))^{-1}=(-\tilde x,h_{\tilde x}(v-\tilde x)^{-1}).$$
Then if $g_i=(\tilde x_i,h_{\tilde x_i}(v))$ for $i=1,2$
$$
\begin{array}{rcl}
[g_1,g_2]&=&g_1g_2g_1^{-1}g_2^{-1} \\
&=&( x_1,h_{\tilde x_1}(v))( x_2,h_{\tilde x_2}(v))
(- x_1,h_{\tilde x_1}(v-\tilde x_1)^{-1})(- x_2,h_{\tilde x_2}(v-\tilde x_2)^{-1})\\
&=&( x_1+ x_2,h_{\tilde x_1}(v)h_{\tilde x_2}(v+\tilde x_1))
(- x_1- x_2,h_{\tilde x_1}(v-\tilde x_1)^{-1}h_{\tilde x_2}(v-\tilde x_1-\tilde x_2)^{-1})\\
&=&(0,h_{\tilde x_1}(v)h_{\tilde x_1}(v+\tilde x_2)^{-1}
h_{\tilde x_2}(v+\tilde x_1)h_{\tilde x_2}(v)^{-1})\\
&=&(0,H_v(\tilde x_1,\tilde x_2)H_v(\tilde x_2,\tilde x_1)^{-1}).\end{array}
$$
Hence $e^{\mathcal{L}}(x_1,x_2):=\iota^{-1}([g_1,g_2])=H_v(\tilde x_1,\tilde x_2)H_v(\tilde x_2,\tilde x_1)^{-1}$.
Note that since the left hand side is independent of $v$, the right hand side is.
\end{proof}

\subsection{Proof of Theorem \ref{dichotomy}}\label{finite}

In the theory of complex tori, given a line bundle $\mathcal{M}$ over a torus $X_\Lambda$ it is shown that a certain group $K(\mathcal{M})$ (analogous to the group we have defined for quantum tori) is either finite, or the whole of $\mathcal{M}$.  The proof of this result relies on the fact that torus $X_\Lambda$ can be viewed as a complete projective variety, and that the corresponding pairing $e^\mathcal{M}$ is a morphism of projective varieties.\newline

Theorem \ref{dichotomy} proves this result for quantum tori.

\begin{proof}[of part \ref{nontriviald}]  
Suppose that $Ch(\mathcal{L})=\eta$ is nontrivial.\\

I first claim that it suffices to only consider those line bundles represented by the cocycles $\hat{\sigma}(\eta)_l(v)$ defined in \eqref{sigmasp}.  If $\mathcal{L}_1$ and $\mathcal{L}_2$ are isomorphic line bundles then it is clear that $K(\mathcal{L}_1)=K(\mathcal{L}_2)$.  \\

For $c \in \B{C}^\times$ we let $\gamma_c$ denote the homomorphism of $L$ defined by 
\begin{equation}\label{gamma}
\gamma_c(a\omega_1+b\omega_2)=c^b.
\end{equation}
By Theorem \ref{theorysplit} it suffices to only consider those line bundles represented by cocycles of the form $\gamma_c\hat{\sigma}_l(v)$ for some $\eta \in \textrm{Alt}^2(L,\B{Z})$ and $c \in \B{C}^\times$.  If $A_l^\mathcal{L}(v)$ is a cocycle of this form then
$$\frac{A_l^\mathcal{L}(v+\tilde x)}{A_l^\mathcal{L}(v)}=\frac{\hat \sigma_l(v+\tilde x)}{\hat \sigma_l(v)}.$$
Therefore if $\mathcal{M}$ corresponds to the cocycle $\hat \sigma_l(v)$ then $K(\mathcal{L})=K(\mathcal{M})$. \newline

We now show that $K(\mathcal{L})$ is isomorphic to $(\B{Z}/s_\eta \B{Z}) \times(\B{Z}/s_\eta \B{Z})$.  If $x=\pi(\tilde x) \in K(\mathcal{L})$ then there exists a unique $h_{\tilde x} \in \mathcal{H}^\ast$ such that $h_{\tilde x}(0)=1$ and 
$$\frac{\hat{\sigma}_l(v+\tilde x) }{ \hat{\sigma}_l(v) }=\frac{h_{\tilde x}(v+l) }{ h_{\tilde x}(v)}.$$
However, explicit calculation shows that
\begin{equation}\label{yelp}
\frac{\hat{\sigma}_l(v+\tilde x) }{ \hat{\sigma}_l(v) }=e^{2 s_\eta\frac{\pi i }{ \omega_1} b \tilde x}
\end{equation}
where $l=a \omega_1+b \omega_2 \in L$. This latter expression is independent of $v$.\newline

There exists a holomorphic function $H_{\tilde x}$ such that $h_{\tilde x}(v)=e^{2 \pi i H_{\tilde x}(v)}$. By \eqref{yelp}, for all $l \in L$, $H_{\tilde x}(v+l)-H_{\tilde x}(v)$ is a holomorphic function independent of $v$, and by continuity is constant on lines of constant imaginary part.  Calculating the derivative of $H_{\tilde x}$ along a line of constant imaginary part we see that $H'(v)=0$. Hence there exist constants $k(\tilde x)$ and $c(\tilde x)$ such that
$$H_{\tilde x}(v)=\frac{k(\tilde x) }{ \omega_1}v+c(\tilde x).$$
Now we compute
$$\frac{h_{\tilde x}(v+l) }{ h_{\tilde x}(v)}=e^{\frac{2 \pi i }{ \omega_1}k(\tilde x)(a\omega_1+b \omega_2)}.$$
Equating this last expression with that of \eqref{yelp} we obtain
$$e^{2{\pi i \over \omega_1}s_\eta b \tilde x} =e^{2 {\pi i \over \omega_1} k(\tilde x)(a \omega_1+b \omega_2)}.$$
Note that the right hand side is dependent on $a$, whereas the left hand side is not.  Since this equality holds for all $a \in \B{Z}$ we therefore have $k(\tilde x) \in \B{Z}$.
We deduce that
$$s_\eta b \tilde x \in \B{Z}(a \omega_1+b \omega_2)+\B{Z} \omega_1. $$
This holds for all $a,b \in \B{Z}$, and hence $s_\eta \tilde x \in L$. \newline

Conversely, if $s_\eta \tilde x \in L$, then we have 
$$\tilde x={\alpha \over s_\eta}\omega_1+{\beta \over s_\eta}\omega_2$$
for some $\alpha,\beta \in \B{Z}$.  Define
$$h_{\tilde x}(v)=e^{{2 \pi i \over \omega_1}\beta v}.$$
Then 
$${\hat \sigma_l(v+\tilde x) \over \hat \sigma_l(v)}={h_{\tilde x}(v+\tilde x) \over h_{\tilde x}(v)}.$$
Hence $\Lambda(\mathcal{L}) \cong {1 \over s_\eta}L$, and the result follows. \newline

\noindent
The explicit formula for $Ch(\mathcal{L})$ follows from Proposition \ref{c2}.\\

\noindent
{\it Proof of part \ref{triviald}.}  Suppose $\mathcal{L}$ is a line bundle over $Z_L$ represented by a cocycle $A_l(v)$.  \\

$\ref{p} \Rightarrow \ref{q}$: Suppose $Ch(\mathcal{L})=0$. By Proposition \ref{constant} $A_l(v)$ is cohomologous to a constant cocycle $K_l(v)$, and hence we can use this representative of the cohomology class to determine $K(\mathcal{L})$:
$$K(\mathcal{L})=\left\lbrace x \in Z_L : (\exists \tilde x \in \pi^{-1}(x)) \land \left( \frac{K_l(v+\tilde x)}{K_l(v)} \in
H^1(L,\mathcal{H}^\ast)\right) \right\rbrace. $$
But since $K_l(v)$ is constant in $v$ we have $K(\mathcal{L})=Z_L$.\newline 

$\ref{q} \Rightarrow \ref{r}$: Suppose $K(\mathcal{L})=Z_L$.  By part \ref{nontriviald} of the theorem, if $\eta=Ch(\mathcal{L})$ is non zero then $\Lambda(\mathcal{L})={1 \over s_\eta}L$. Since $K(\mathcal{L})=Z_L$ we have $\Lambda(\mathcal{L})=\B{R}$, so we must have $\eta=0$.  Hence $\mathcal{L}$ is given by a cocycle of the form
$$A_l(v)=\gamma_c(l){h(v+l) \over h(v)}$$
for some $h \in \mathcal{H}^\ast$.  We see that for $\tilde x \in \B{R}$
$${A_l(v+\tilde x) \over A_l(v)}={h(v+\tilde x+l)\over h(v+l)} {\Bigg \slash} {h(v+\tilde x) \over h(v)}.$$
By the uniqueness of $h_{\tilde x}(v)$ we find
$$h_{\tilde x}(v)={h(v+\tilde x) \over h(v)}{h(0) \over h(\tilde x)}.$$
Using the explicit formula for $H_v(\tilde x_1,\tilde x_2)$ in Definition \ref{lambdadefine} we have
$$H_v(\tilde x_1,\tilde x_2)={h(v+\tilde x_1+\tilde x_2) \over h(v+\tilde x_1)h(v+\tilde x_2)h(v)}.$$
Observe that this is symmetric in $\tilde x_1$, and $\tilde x_2$, and hence by Proposition \ref{c2} we conclude that $e^{\mathcal{L}}(\tilde x_1,\tilde x_2)=1$. 
\newline

$\ref{r} \Rightarrow \ref{p}$:  Suppose $e^\mathcal{L} \equiv 1$, and assume that $Ch(\mathcal{L}) \neq 0$.  Then $s_\eta \neq 0$, and by the results of part \ref{nontriviald} we have $Ch(\mathcal{L})(\omega_1/s_\eta,\omega_2/s_\eta)=1/s_\eta$, and hence
$e^{\mathcal{L}}(\omega_1/s_\eta,\omega_2/s_\eta)=e^{-2 \pi i/s_\eta}\neq 1$. This is a contradiction, and hence we must have $Ch(\mathcal{L})=0$. 
\end{proof}
As an immediate corollary we have

\begin{cor}
Let $\mathcal{L}$ be a line bundle over $Z_L$ with Chern class $\eta$.
Extend $Ch(\mathcal{L})$ to a pairing on ${1 \over s_\eta}L$ by linearity, and let
$\tilde x_1,\tilde x_2 \in \Lambda(\mathcal{L})$.  Then
$$e^\mathcal{L}(\tilde x_1,\tilde x_2)=e^{2 \pi i Ch(\mathcal{L})(\tilde x_1,\tilde x_2)}.$$
\end{cor}
\begin{proof}
By Proposition \ref{c2} we have
$$e^\mathcal{L}(\tilde x_1,\tilde x_2)={h_{\tilde x_2}(v+\tilde x_1) \over h_{\tilde x_2}(v)}
{h_{\tilde x_1}(v) \over h_{\tilde x_1}(v+\tilde x_2)}.$$
Let $s_\eta\tilde x_1=a\omega_1+b\omega_2$ and $s_\eta \tilde x_2=c\omega_1+d \omega_2$.  Then
$$
\begin{array}{rcl}
e^\mathcal{L}(\tilde x_1,\tilde x_2)&=&e^{2{ \pi i \over s_\eta \omega_1} d(a\omega_1+b\omega_2)}e^{-2{ \pi i \over s_\eta \omega_1} b(c\omega_1+d\omega_2)}\\
&=&e^{2 \pi i{ad-bc \over s_\eta}}\\
&=&e^{2 \pi i Ch(\mathcal{L})(\tilde x_1,\tilde x_2)}.
\end{array}
$$
\end{proof}

\subsection{Theta functions and Hilbert's twelfth problem for real quadratic fields}

Having defined a notion of line bundles for quantum tori, it is natural to question the existence of sections of such line bundles.  Given a line bundle $\mathcal{L}$ given by $A_l^\mathcal{L}(v) \in Z^1(L,\mathcal{H}^\ast)$, a holomorphic section of $\mathcal{L}$ is given by a homorphic function $\theta(v)$ such that for all $l \in L$
\begin{equation}\label{they}
\theta(l+v)=A_l^\mathcal{L}(v)\theta(v).
\end{equation}
A function satisfying this property will be called a \emph{theta function} for $Z_L$.  In the analogous theory for complex tori, the Jacobi theata functions satisfy a corresponding functional equation, which play a central role in the solution of Hilbert's twelfth problem for imaginary quadratic fields. 

\begin{prop}
Let $\theta$ be a nonconstant holomorphic theta function for a quantum torus $Z_L$.  Then there exists $h \in \mathcal{H}^\ast$ and constants $A \in \mathbb{C}^\times$, $\alpha \in \B{C}$ such that
$$\theta(v)=h(v)e^{2 \pi i \alpha v}.$$
\end{prop}
\begin{proof}
First we note that $\theta$ can have no zeros, for otherwise by \eqref{they} it would have a dense set of zeros along a line in $\B{C}$ and hence be identically zero.\\

Now suppose $\theta$ satisfies \eqref{they} and $$B_l(v)=A_l^\mathcal{L}(v)\frac{h(v+l)}{h(v)}$$ for some non-vanishing holomorphic function $h$, then
$\theta(v)h(v)$ is a theta function for the line bundles associated to $B_l$.  Hence it suffices to show that there are no nonconstant holomorphic theta functions satisfying \eqref{they} for a representative of each cohomology class in $Z^1(L,\mathcal{H}^\ast)$.  By Theorem \ref{theorysplit} we need only consider line bundles represented by cocycles of the form
$$\gamma_c\hat \sigma(\eta)_l(v)$$
as defined in equations \eqref{sigmasp} and \eqref{gamma}.\\

Note that if $\theta$ is a theta function for this cocycle then we have
\begin{equation}\label{lim}
\left|\theta(v+l) \right|=\left|c^b\hat \sigma(\eta)_l(v) \theta(v)\right|=\left| c \right|^b\left| \theta(v)\right|.
\end{equation}
We may choose a sequence $l_n=a_n\omega_1+b_n\omega_2 \in L$ such that $l_n \rightarrow 0$ but $b_n \rightarrow \infty$ as $n \rightarrow \infty$.  If $\left| c \right| \neq 1$ then \eqref{lim} gives us a contradiction.  Hence $\left|c \right|=1$, and for all $v \in \B{C}$, $l \in L$ we have 
\begin{equation}\label{periodl}
\left|\theta(v+l) \right|=\left| \theta(v)\right|.
\end{equation}

Now fix $r \in \B{R}$.  Since $\theta$ is nonvanishing there exists a function $x_r(v)$ holomorphic in $v$ such that
\begin{equation}\label{s}
{\theta(v+r) \over \theta(v)}=e^{2 \pi i x_r(v)}.
\end{equation} 
Without loss of generality we assume that $x_0(v)=0$.  Equation \eqref{periodl} implies that $x_r(v) \in \B{R}$ for all $v \in \B{C}$.  Since $x_r(v)$ is a holomorphic function in $v$ this implies that it is constant. \newline

Now fix $v \in \B{C}$, and let $r,s \in \B{R}$.  
Then 
$${\theta(v+r+s) \over \theta(v)}={\theta(v+r+s) \over \theta(v+r)}{\theta(v+r) \over \theta(v)}. $$
Hence there exists $n(v) \in \B{Z}$ such that
$$
\begin{array}{rcl}
x_{r+s}(v)&=&x_s(v+r)+x_r(v)+2 \pi i n(v)\\
&=&x_s(v)+x_r(v) + 2 \pi i n(v).
\end{array}
$$
Since $x_0(v)=0$ we see that $n(v)=0$, and as a function of $r \in \B{R}$, $x_r(v)$ is a homomorphism. Hence for all $r \in \B{R}$, $x_r(v)=\alpha r$ for some $\alpha \in \B{R}$.\newline

Now consider the left hand side of \eqref{s}.  As $r$ varies over $\B{C}$ this is a holomorphic function.  Hence 
for fixed $v$, there exists a function $x_v(w)$ holomorphic in $w$ such that for all $v \in \B{C}$
\begin{equation}\label{label}
{\theta(v+w) \over \theta(v)}=e^{2 \pi i x_v(w)}.
\end{equation}
Again we assume without loss of generality that $x_v(0)=0$, and therefore $x_z(v)=x_v(z)$ for all $v, w \in \B{C}$. On $\B{R}$ we therefore have $x_v(w)=\alpha w$, and hence by holomorphicity holds on the whole plane.  Putting $A=\theta(0)$ we see from \eqref{label} that $\theta(v)=Ae^{2 \pi i \alpha v}$.
\end{proof}

By this result up to a nonvanishing holomorphic function the holomorphic theta functions for quantum tori are simply exponential functions, which seem unlikely to have a central role in generating class fields of real quadratic fields.  \\

However, in this paper we have considered line bundles defined by elements of $Z^1(L,\mathcal{H}^\ast)$.  Let $\mathcal{K}$ denote the field of fractions of $\mathcal{H}$ and consider the group $Z^1(L,\mathcal{K}^\ast)$.  In \cite{ShintaniIII}, Shintani considers the \emph{double sine function}, which can be viewed as a theta function $S$ satisfying a relation
$$S(v+l)=A_l(v)S(v)$$
for some $A_l(v) \in Z^1(L,\mathcal{K}^\ast)$.  He uses special values of this function to generate abelian extensions of real quadratic fields in certain special cases.  This example suggests that the notion of holomorphic line bundles is too strong for their application to Hilbert's problem, and that further knowledge of the group $Z^1(L,\mathcal{K}^\ast)$ would be valuable in the context of this problem.

\subsection{Line bundles in Model Theory}

If $X_\Lambda$ is a complex torus, it is the existence of an integral valued alternating form on the lattice $\Lambda$ that ensures the existence of ``very ample'' line bundles over $X_\Lambda$.  Fundamentally, it is the existence of these line bundles that imply that $X_\Lambda$ can be viewed as an algebraic variety.\\

For a line bundle over a complex torus $X_\Lambda=\B{C}/\Lambda$, we have an analogous theory to that developed for quantum tori.  For each line bundle $\mathcal{M}$ over $X_\Lambda$ we obtain a subgroup $K(\mathcal{M})$ of $X_\Lambda$ and an alternating pairing
$$e^\mathcal{M}: K(\mathcal{M}) \times K(\mathcal{M}) \rightarrow \B{C}^\ast.$$
The group $K(\mathcal{M})$ is an algebraic subvariety of $X_\Lambda$, and the pairing $e^\mathcal{M}$ is a morphism of algebraic varieties.\\

Although quantum tori are not algebraic varieties, through the work of Zilber they can be defined in a class of structures called \emph{Analytic-Zariski structures}.  These structures represent a variation to the notion of a Zariski structure introduced by Zilber and Hrushovski in \cite{Hrushovski}.
\newline

The structure considered by Zilber is the two sorted structure 
\begin{equation}\label{structure}
{T}_\theta:=((\B{C},+,{A}_\theta),\exp,\B{C}^\ast)
\end{equation}
where
$${A}_\theta=(i+\theta)\B{R}+ 2 \pi \B{Z}+2 \pi i \B{Z}.$$
It is shown in \cite{ZilberII} that that the quotient $\B{C}/{A}_\theta$ is isomorphic to a quantum torus.  Letting $\mathcal{G}_\theta=\exp(A_\theta)$, and considering the structure $P_\theta:=(\B{C}^\times, \mathcal{G}_\theta, .)$, the theory of this structure is known to be superstable providing Schanuel's conjecture holds. The quantum torus associated to the pseudolattice $L_\theta:=\B{Z}+\B{Z}\theta$ can be interpreted in this structure as the quotient $\B{C}^\times/\mathcal{G}_0$. \\

\begin{prop}\label{az}
If Shanuel's conjecture holds, then $K(\mathcal{L})$ is an Analytic-Zariski set.
\end{prop}
\begin{proof}
We consider the cases $s\neq 0$ and $s =0$ separately.  In the case when $s \neq 0$, by part \ref{nontriviald} of Theorem \ref{dichotomy} we have
$$\Lambda(\mathcal{L}) \cong \{x \in \B{C}^\ast : x^s \in \mathcal{G}_\theta \}.$$
$K(\mathcal{L})$ is the quotient of this by $\mathcal{G}_\theta$, and hence a definable subset of the structure
$(\B{C}^\ast, \mathcal{G}_\theta,.)$.\\

If $s=0$, by part \ref{triviald} of Theorem \ref{dichotomy} we have $K(\mathcal{L}) \cong \mathbf{T}_\theta=\B{C}^\ast/\mathcal{G}_\theta$, which is an Analytic-Zariski structure by Zilber's study (modulo Schanuel's conjecture).  
\end{proof}

At present, this is the limit to the extent we can achieve our goal of defining the pairing $e^\mathcal{L}$ in an Analytic-Zariski structure..  The graph of $e^\mathcal{L}$ is not definable in any of the structures that Zilber considers in \cite{ZilberIII}, \cite{ZilberI} and \cite{ZilberII}.  To achieve this it would be desirable to have a log-function between $\B{C}^\ast$ and $\B{C}$.  We do not know whether the addition of this function to any of Zilber's structures would alter the stability of such a structure.  However, in \cite{ZilberII} Zilber defines a ``random logarithm'' function from $\B{C}^\ast$ to $e^{2 \pi \theta \B{Z}}$ resulting in an unstable structure.

\end{document}